\documentclass{amsart}

\usepackage{amsmath}
\usepackage{amssymb}
\usepackage{amsfonts}

\newtheorem{theorem}{Theorem}
\newtheorem{prop}{Proposition}
\newtheorem{lemma}{Lemma}
\newtheorem{remark}{Remark}
\newtheorem{cor}{Corollary}
\newtheorem{example}{Example}

\newtheorem{definition}{Definition}

\def\Z{\mathbb Z}

\def\Q{\mathbb Q}
\def\C{\mathbb C}
\def\bP{\mathbb P}

\def\L{\mathcal L}

\def\H{\mathcal H}
\def\M{\mathcal M}

\def\l{\lambda}

\def\G{\Gamma}
\def\a{\alpha}
\def\b{\beta}
\def\p{\mathfrak p}
\def\P{\mathcal P}
\def\e{\varepsilon}
\def\iso{\equiv}

\def\lar{\to}
\def\bG{\overline G}
\def\g{\gamma}
\def\bg{\bar \gamma}
\def\u{\mathfrak u}

\def\iso{{\, \cong\, }}

\def\<{\langle}
\def\>{\rangle}

\def\X{\mathcal X}
\def\Aut{\mbox{Aut }}

\def\bAut{\overline{\mbox{Aut}} \, \, \, }

\def\emb{\hookrightarrow }

\def\tr{{\mbox tr }}
\def\det{{\mbox det }}

\def\D{\Delta }
\def\aa{\mathfrak a}

\begin{document}
\title{Hyperelliptic curves with extra involutions}

\date{June  29, 2004 and, in revised form, June  ?? , 2004.}

\author{J. Gutierrez}
\email{jaime@matesco.unican.es}
\address{Department of Mathematics\\ Universidad  de Cantabria\\ E-39071\\ Santander\\ Spain}

\thanks{The first author was supported in part by a grant from the Spain Ministry Science,  PR2002-0009.}

\author{T. Shaska}
\email{tshaska@math.uci.edu}
\address{300 Brink Hall\\ Department of Mathematics\\ University of Idaho\\  Moscow, ID, 83844 \\ USA}
\thanks{The second author was supported in part by UI research office, S0511.}



\begin{abstract}
The purpose of this paper is to study hyperelliptic curves with extra involutions.   The locus $\L_g$ of such
genus $g$ hyperelliptic curves is a $g$-dimensional subvariety of the moduli space of hyperelliptic curves
$\H_g$. We discover a birational parametrization of $\L_g$ via dihedral invariants and show how  these
invariants can be used to determine  the field of moduli of points $\p \in \L_g$.

We conjecture that for $\p\in \H_g$ with $|\Aut(\p)| > 2$ the field of moduli is a field of definition and prove
this conjecture for  any point $\p\in \L_g$ such that the Klein 4-group is embedded in the reduced automorphism
group of $\p$. Further, for  $g=3$ we  show that  for every moduli point $\p \in \H_3$ such that $| \Aut (\p) | >
4$, the field of moduli is a field of definition and provide a rational model of the curve over its field of
moduli.
\end{abstract}
\maketitle
\section{Introduction}

It is an interesting problem in algebraic geometry   to obtain a generalization of the theory of elliptic
modular functions to the case of higher genus. In the elliptic case this is done by the so-called
$j$-\textit{invariant} of elliptic curves. In the case of genus $g=2$, Igusa (1960) gives a complete solution
via {\it absolute invariants} $i_1, i_2, i_3$ of genus 2 curves. Generalizing such results to higher genus is
much more difficult due to the existence of non-hyperelliptic curves. However, even restricted to the
hyperelliptic moduli $\H_g$ the problem is still unsolved for $g \geq 3$. In other words, there is no known
way of identifying isomorphism classes of hyperelliptic curves of genus $g\geq 3$. In terms of classical
invariant theory  this means that the  field of invariants of binary forms of degree $2g+2$ is not known for
$g\geq 3$.

In this paper we focus on the locus $\L_g$ of genus $g$ hyperelliptic curves with extra (non-hyperelliptic)
involutions defined over an algebraically closed field $k$. We determine invariants that generically identify
isomorphism classes of curves in $\L_g$.  Eq.~\eqref{eq} gives a normal form for genus $g$ hyperelliptic
curves with extra involutions. This normal form depends on parameters $a_1, \dots , a_g\in k$. We discover an
action of the dihedral group $D_{g+1}$ of order 2g+2 that symmetrizes $a_1, \dots , a_g$. Invariants of this
action are parameters $ u_1, \dots , u_g\in k [ a_1, \dots , a_g ]$. We call such invariants \emph{dihedral
invariants}  of hyperelliptic curves and show  that $k(\L_g)=k(a_1, \dots , a_g)^{D_{g+1}}$. More precisely,
this $g$-tuple of dihedral invariants parameterizes isomorphism classes of genus $g$ hyperelliptic curves
with extra involutions. The map $k^g\setminus \{\Delta \neq 0\} \to \L_g$ is birational. Thus, dihedral
invariants $u_1, \dots , u_g$ yield a birational parametrization of the locus $\L_g$. Computationally these
invariants give an efficient way  of determining a point of the moduli space $\L_g$. Normally, this is
accomplished by invariants of $GL_2(k)$ acting on the space of binary forms of degree $2g+2$. These
$GL_2(k)$-invariants are not known  for $g \geq 3$.  However, dihedral invariants are explicitly defined  for
all $g$. The most direct method to compute the dihedral invariants  requires the curve in the normal form.
This can be done by solving a polynomial system of equations.

In section 4, we study the field of moduli of hyperelliptic curves in $\L_g$.  Whether or not the field of
moduli is a field of definition is in general a difficult problem that goes back  to Weil, Baily, Shimura et
al. We conjecture that for each $\p \in \H_g$ such that $|\Aut (\p)| > 2$ the field of moduli is a field of
definition. Again we focus only on the locus $\L_g$. Making use of $(\u_1, \dots , \u_g)$,  we show that if
the Klein 4-group can be embedded in the reduced automorphism group of $\p\in \L_g$ then the conjecture
holds. Moreover, the field of moduli is a field of definition for all $\p \in \L_3$ such that $|\Aut
(\p)|>4$.

\medskip

\noindent   \textbf{Notation:}  Throughout this paper $k$ denotes an algebraically closed field of
characteristic not equal to 2,   $V_4$ denotes the Klein 4-group,  $D_{n}$ (resp., $\Z_n$) the  dihedral
group of order $2n$ (resp., cyclic group of order $n$), and $\Gamma:=PGL_2(k)$.
%
\section{Preliminaries}
Let $k(X)$ be the field of rational functions in $X$. We identify the places of $k(X)$ with the points of
$\bP^1=k\cup\{\infty\}$ in the natural way (the place $X=\a$ gets identified with the point $\a\in \bP^1$).
Let $K$ be  a quadratic extension  field of $k(X)$  ramified exactly  at $n$ places $\a_1, \dots , \a_n$ of
$k(X)$. The corresponding places of $K$ are called the  \textbf{Weierstrass points}  of $K$.  Let $\P:=\{
\a_1, \dots , \a_n \}$. Then $K=k(X,Y)$, where
\begin{equation}\label{2.0}
Y^2=\prod_{\overset {\a\in \P}{ \a \neq \infty} }  (X-\a).
\end{equation}
Let $G=Aut(K/k)$. It is well known that $k(X)$ is the only  genus 0 subfield of degree 2 of $K$; thus $G$
fixes $k(X)$. Thus, $G_0:=Gal(K/k(X))= \< z_0 \>$, with $z_0^2=1$,  is  central in $G$. We call the
\textbf{reduced automorphism group}  of $K$ the group $\bG:=G/G_0$.  Then, $\bG$ is naturally isomorphic to
the subgroup of $Aut(k(X)/k)$ induced by $G$. We have a natural  isomorphism $\Gamma:=PGL_2(k) {\overset \iso
\lar} Aut(k(X)/k)$.

The action of $\Gamma$ on the places of $k(X)$ corresponds under the above identification  to the usual action on
$\bP^1$ by fractional linear transformations: $t \mapsto \frac {at+b} {ct+d}$. If $l$ is prime to $char(k)$ then
each  element of order $l$ of $\Gamma$ is conjugate to $\begin{pmatrix} {\e_l} & 0 \\ 0 & 1 \end{pmatrix}$, where
$\e_l$ is a  primitive $l$-th root of unity.  Each such element has 2 fixed points on $\bP^1$ and other orbits are
of length $l$.  If $l=char(k)$ then, $\Gamma$ has exactly one class of elements of order $l$, represented by
$\begin{pmatrix} 1 & 1\\ 0 & 1\end{pmatrix}$. Each such element has exactly one fixed point on $\bP^1$. Further,
$G$ permutes $\a_1, \dots , \a_n$.  This yields an embedding $\bG \hookrightarrow S_n$.
\begin{lemma}\label{lem0}
Let $\g \in G$ and $\bg$ its image in $\bG$.  Suppose $\bg$ is an involution. Then, $\g$ has order 2 if and only
if it fixes no Weierstrass points.
\end{lemma}
\begin{proof} Suppose $\bg $ is an involution. By the above we may assume $\bg(X)=-X$. We may further assume
that $1\in \P$ by replacing $X$ by $cX$ for a suitable $c\in k^*$. Now  assume $\bg$ fixes no points in $\P$.
Thus,
%
$\P=\{ \pm 1, \pm \a_1,    \dots , \pm \a_{\frac {n-2} 2}  \},$
%
where $\a_i \in \bP^1\setminus \{0,\infty, \pm 1\}$. Hence,
%
$ Y^2=(X^2-1) \prod_{i=1}^{g} (X^2-\a_i^2)$,   for $ g= \frac {n-2} 2.$
So, we have $\g(Y)^2=Y^2$. Hence $\g(Y)=\pm \, Y$, and $\g$ has order 2.

Suppose  $\bg$  fixes 2 points of $\P$.  Then, $\P=\{0,\infty, \pm 1, \pm \a_1,  \dots , \pm \a_s  \}$, where
$\a_i \in \bP^1\setminus \{0, \infty, \pm 1\}$.  Hence,
$Y^2=X(X^2-1) \prod_{i=1}^{g} (X^2-\a_i^2)$,  for $g= \frac {n-4} 2$.
So $\g(Y)^2=- Y^2 $ and $\g(Y)=\sqrt{-1} \,\, Y$. Hence, $\g$ has order 4.
\end{proof}

Because $K$ is the unique degree 2 extension of $k(X)$ ramified  exactly at $\a_1$, \dots  , $\a_n$, each
automorphism of $k(X)$ permuting these $n$ places extends to  an automorphism of $K$. Thus, $\bG$ is the
stabilizer in $Aut(k(X)/k)$ of the set $\P$. Hence under the isomorphism   $\Gamma \mapsto  Aut(k(X)/k)$, $\bG$
corresponds to the stabilizer $\Gamma_\P$ in $\Gamma$ of the $n$-set $\P$.

An  \textbf{extra involution} of $K$ is an involution in $G$ which is different from $z_0$ (the hyperelliptic
involution). If $z_1$ is an extra involution and $z_0$ the hyperelliptic one, then $z_2:=z_0\, z_1$  is
another extra involution in $G$. So the extra involutions come naturally in pairs.  These pairs correspond
bijectively to pairs $F_1, F_2$ of degree 2 subfields of $K$ with $F_1 \cap k(X) = F_2 \cap k(X)$.  An
involution in $\bG$ is called extra involution if it is the image of an \textbf{extra involution of $G$}.
\begin{lemma} \label{lem1}
Suppose $z_1$ is an extra involution of $K$. Let $z_2 :=z_1\, z_0$, where $z_0$ is the hyperelliptic involution.
Let $F_i$ be the fixed field of $z_i$ for $i=1,2$. Then $K=k(X,Y)$ where
\begin{equation}\label{eq}
Y^2=X^{2g+2} + a_{g} X^{2g}+ \dots + a_1 X^2 +1
\end{equation}
and  $\Delta (a_1, \dots , a_g) \neq 0$ (i.e., $\Delta$ is the discriminant of the right hand side).
Furthermore,  $F_1$ and $F_2$ are the subfields $k(X^2,Y)$ and $k(X^2, YX)$.
\end{lemma}
\proof Recall that $z_0(X, Y)=(X, -Y)$. We choose the coordinate $X$ such that  $\bar z_1(X) =-X$. By Lemma
1,  the involution $z_1$ fixes no points of $\P$, hence $\P=\{\pm \a_1, \dots  , \pm \a_s \}$, where $s= g+1$
and $\a_i \in k\setminus \{0\}$.

Let $\b_i := \a_i^2$, for $i=1, \dots , s$. Then
we have $K=k(X,Y)$ with $Y^2=\prod_{i=1}^s (X^2-\b_i).$
Let $a_1, \dots , a_g$ denote symmetric polynomials of $\b_i$ (up to a sign change). Then,
\begin{equation}
Y^2=X^{2g+2} + a_{g} X^{2g}+ \dots + a_1 X^2 +a_0.
\end{equation}
We may further replace $X$ by $\l X$, for a suitable $\l$, to get $a_0=(-1)^s \prod_{i=1}^s \b_i=1$.

Since the roots $\a_1, \dots, \a_s$ are distinct then $\Delta (a_1, \dots , \a_g) \neq 0$ (i.e., $\Delta$ is
the discriminant of the right hand side). The elements $X^2$ and $XY$ are fixed by $z_2$. This implies the
claim. \qed

We will consider pairs $(K, z)$ with $K$ a genus $g$ field and $z$ an extra involution.  Two such pairs $(K,
z)$ and $(K', z')$ are called \textbf{isomorphic} if there is a $k$-isomorphism $\a : K \to K'$ with $z'= \a
z \a^{-1}$. Determining these isomorphism classes will be the focus of the next section.

\section{Dihedral invariants}
Let $\X_g$  be  a hyperelliptic curve of genus $g\geq 2$ defined over $k$ and $K$ its function field. Then, $\X_g$
can be described as a double cover of $\bP^1:=\bP^1(k)$ ramified in $( 2g+2)$ places $w_1, \dots , w_{2g+2}$. This
sets up a bijection between isomorphism classes of hyperelliptic genus $g$ curves and unordered distinct
$(2g+2)$-tuples $w_1, \dots , w_{2g+2} \in \bP^1 $ modulo automorphisms of $\bP^1 $. An unordered $(2g+2)$-tuple
$\{w_i\}_{i=1}^{2g+2}$ can be described by a binary form (i.e. a homogeneous equation $f(X,Z)$) of degree
$(2g+2)$. Hence, we assume that $\X_g$ is  given by
%
$Y^2 Z^{2g}=f(X,Z)=\sum_{i=0}^{2g+2} a_i X^i Z^{2g+2-i}.$

Let $\H_g$ denote the moduli space of hyperelliptic genus $g$ curves.  To describe $\H_g$ we need to find
rational functions of the coefficients of a binary form $f(X,Z)$, invariant under linear substitutions in
$X,Z$. Such functions are traditionally called  \textit{absolute invariants} for $g=2$; see Igusa \cite{Ig}
or Krishnamorthy/Shaska/V\"olklein \cite{vishi}. We will adapt the same terminology even for $g \geq 3$. The
absolute invariants are $GL_2(k) $ invariants under the natural action of $GL_2(k)$ on the space of binary
forms of degree $2g+2$. Two genus $g$ hyperelliptic curves are isomorphic if and only if they have the same
absolute invariants. We denote by $\L_g$ the locus in $\H_g$ of hyperelliptic curves with extra involutions.
To find an explicit description of $\L_g$  means finding explicit equations in terms of absolute invariants.
Such equations are computed only for $g=2$; see Shaska/V\"olklein \cite{Sh-V}.  Computing similar equations
for $g\geq 3$ requires first finding the corresponding absolute invariants. This is still an open problem in
classical invariant theory even for $g=3$. Even in the case that absolute invariants are known, they are
expected to have very large expressions in terms of coefficients of the binary forms. Thus, equations
defining $\L_g$ are expected to be very large and not helpful for any practical use. In this section we find
new parameters for $\L_g$, which we call dihedral invariants. This $g$-tuple  $\u \in k^g$ generically
classifies isomorphism classes of curves $\X_g \in \L_g$.
\subsection{The dihedral group action on $k(a_1, \dots a_g)$}
Let $\X_g$ be a genus $g$ hyperelliptic curve with an extra involution. Then,  $\X_g$ is given as  in
Eq.~\eqref{eq}. We need to determine to what extent the normalization  in the  proof of Lemma 2 determines
the coordinate $X$.

The condition $z_1(X)=-X$ determines the coordinate $X$ up to a coordinate change by some $\gamma \in \Gamma$
centralizing $z_1$. Such $\gamma$ satisfies $\gamma (X)=mX$ or $\gamma (X) = \frac m X$, $m \in k\setminus \{0\}$.
The additional condition $(-1)^g \b_1 \cdot \cdot \cdot \cdot \b_{g+1}=1$ forces
\begin{equation}
(-1)^g \, \g(\a_1) \dots \g(\a_{2g+2})=1.
\end{equation}
Hence,   $m^{2g+2}=1$. So $X$ is determined up to a coordinate change by the subgroup $ D_{g+1} < \Gamma$
generated by $\tau_1: X\to \e X$, $\tau_2: X\to \frac  1 X$, where $\e$ is a primitive $(2g+2)$-th root of
unity. Hence, $D_{g+1}$ acts on $k(a_1, \dots , a_g)$ as follows:
\begin{equation}
\begin{split}
\tau_1: \quad  & a_i \to \e^{2i} a_i, \qquad  for \quad i=1, \dots , g \\
\tau_2: \quad  & a_i \to a_{g+1-i}, \qquad for \quad  i=1, \dots , [\frac {g+1} 2].\\
\end{split}
\end{equation}
Thus, the fixed field $k(a_1, \dots , a_g)^{D_{g+1}}$ is the same as the function field of the variety
$\L_g$. We summarize in  the following:
\begin{prop}
For a fixed genus $g\geq 2$, let $\L_g$ denote the locus of genus $g$ hyperelliptic curves with extra
involutions. Then,  $k(\L_g)=k(a_1, \dots , a_g)^{D_{g+1}}$.
\end{prop}
Next we find the invariants of such action explicitly. The proof of the following lemma is obvious.
\begin{lemma}
Fix  $g\geq 2$.  The following
\begin{equation}
u_i:=  a_1^{g-i+1} \, a_i \, + \, a_g^{g-i+1} \, a_{g-i+1}, \quad for \quad 1 \leq i \leq g\\
\end{equation}
are invariants under the $D_{g+1}$-action and are called   \textbf{dihedral invariants}  of the genus  $g$.
\end{lemma}
It is easily seen that $\u:=(u_1, \dots , u_g)=(0, \dots , 0)$ if and only if $a_1 = a_g = 0$. In this case
replacing $a_1, a_g$ by $a_2, a_{g-1}$ in the formula above would give new invariants. For the rest of the
paper we will focus in the case that $\u \neq 0$, as the other cases are simpler.
For small $g$ (i.e., $g=2,3$), we have the following.
\begin{example} For genus 2, the  dihedral invariants are
\begin{equation}u_1=a_1^3 +a_2^3, \quad   \quad  u_2=2 a_1 a_2, \end{equation}
see \cite{Sh-V} for a detailed study of this case.
For $g=3$  we have
\begin{equation}
u_1 =a_1^4+ a_3^4, \quad   u_2=(a_1^2 + a_3^2)\,  a_2, \quad u_3=2 a_1 a_3.
\end{equation}
\end{example}

The next theorem shows that the dihedral invariants generate $k(\L_g)$, therefore $\L_g$ is a rational variety.
\begin{theorem}
Let $g\geq 2$ and  $\u=(u_1, \dots , u_g)$ be the  $g$-tuple of dihedral invariants. Then, $k(\L_g)=k(u_1, \dots ,
u_g)$.
\end{theorem}
\proof The dihedral invariants are fixed by the $D_{g+1}$-action. Hence, $k(\u)\subset k(\L_g) $. Thus, it is
enough to show that $[k(a_1, \dots a_g): k(\u)]=2g+2$. For each $2 \leq i \leq g-1$ we have
\begin{equation}
\begin{split}
a_1^{g+1-i} a_i + a_g^{g+1-i} a_{g+1-i} = u_i\\ a_1^i a_{g+1-i} + a_g^i a_i=u_{g+1-i}
\end{split}
\end{equation}
giving $\, a_i, \, \, a_{g+1-i} \in k(\u, a_1, a_g)$. Then, the extension $k(a_1, \dots , a_g) /  k(u_1, \dots ,
u_g)$ has equation
\begin{equation}\label{disc}
2^{g+1}\, a_g^{2g+2} - 2^{g+1}\, u_1 \,  a_g^{g+1} + u_g^{g+1}=0
\end{equation}
This  completes the proof.

\qed

\noindent The map
$$\theta: \, \, (a_1, \dots , a_g) \longrightarrow (u_1, \dots , u_g)$$
is a branched Galois covering with group $D_{g+1}$ of the set
$$\{ (u_1, \dots , u_g) \in k^g : \D_\u \neq 0\}$$
by the corresponding open subset of $(a_1, \dots , a_g)$-space, where $\D_\u$ is the discriminant in Lemma~2,
in terms of the dihedral invariants. If $(a_1, \dots , a_g)$ and  $(a_1^\prime, \dots , a_g^\prime)$ have the
same invariants $(u_1, \dots , u_g)$  then they are    $D_{g+1}$ conjugate.

\begin{lemma}
If $\, \aa := (a_1, \dots , a_g)\in k^g$ with $\D_\aa \neq 0$ then Eq.~\eqref{eq} defines a genus $g$ field
$K:=k(X, Y)$ such that its reduced automorphism group contains the extra involution $z_1: X \to -X$. Two such
pairs $(K_\aa, z_1)$ and  $(K_{\aa^\prime}, z_1^\prime)$ are isomorphic if and only if
the corresponding dihedral invariants are the same.
\end{lemma}

\proof The first part of the lemma is obvious as it is the existence of the extra involution $z_1: X \to -X$.
If two pairs are isomorphic then there is $\a : K_\aa \to K_{\aa^\prime}$ which yields $K=k(X, Y)=k(X', Y')$
with $k(X)=k(X')$ such that $X, Y$     satisfy Eq.~\eqref{eq} and $X', Y'$ satisfy the corresponding equation
with $a_1, \dots , a_g$ replaced by $a_1^\prime, \dots , a_g^\prime$. Furthermore, $z_1 (X')= - X'$. Hence
$X'$ is conjugate to $X$ under $\< \tau_1, \tau_2\>$. Thus the dihedral invariants are the same since they
are fixed by $\< \tau_1, \tau_2\>$. The converse goes similarly.
\endproof

\noindent The following theorem is an immediate consequence of the above lemma.

\begin{theorem}
The tuples $\u=(u_1, \dots , u_g)\in k^g$ with $\D\neq 0$ bijectively classify the isomorphism classes of
pairs $(K, z)$ where $K= k(\X_g)$ and $z$ is an involution in $\bAut (\X_g)$. In particular, a given curve
will have as many tuples of these invariants as its reduced automorphism  group has conjugacy classes of
extra involutions.
\end{theorem}

For hyperelliptic curves of genus $g=3, 4$ all tuples of invariants and their algebraic relations
are determined.  For curves with automorphism group isomorphic to $V_4$ we have the following:

\begin{cor}
Let $\X_g$ and $\X_g^\prime$ be genus $g$ hyperelliptic curves with automorphism group isomorphic to $V_4$.
Then, $\X_g$ is isomorphic to $\X_g^\prime$ if and only if they have the same dihedral invariants.
\end{cor}

\proof  Immediate consequence of the above theorem since in this case the reduced automorphism group is
$\Z_2$.

\qed

In general, the case  where the reduced automorphism group  has more involutions can be
characterized in the following:

\begin{theorem} Let $\X_g $ be a genus $g$ hyperelliptic curve with an extra involution and
$(u_1, \dots , u_g)$ its corresponding dihedral invariants.  \\

i) If $V_4\emb     \bAut (\X_g)$ then
$ 2^{g-1}\, u_1^2 = u_g^{g+1}$.

ii) Moreover, if $g$ is odd then $V_4 \emb \bAut( \X_g$ implies that
$$\left( 2^r \, u_1 - u_g^{r+1}\right) \, \left( 2^r \, u_1 + u_g^{r+1}\right)=0$$
where $r= \left[ \frac {g-1} 2 \right]$. The first factor corresponds to the case when  involutions
of $V_4\emb \bar G$ lift to involutions in $G$, the second factor corresponds to the case when two
of the involutions of $V_4\emb \bar G$ lift  to elements of order 4 in $G$.
\end{theorem}

\begin{proof} Since $\X_g$ has an extra involution then it has an equation  as in Eq.~\eqref{eq}.
Moreover, this extra involution in $\bar G$ is given by $ z_1 (X) = -X$ and fixes no Weierstrass points of
$\X_g$; see the proof of Lemma 1.

Let $V_4 \emb \bar G = \bAut (\X_g)$. Then there is  another involution $z_2\neq z_1$   in $\bar G$ such that
$V_4=\<z_1, z_2\>$. Let $M \in \G$ be the corresponding matrix for $z_2$. Then $\tr\, (M) = 0$ and $\det \,
(M) = -1$. Since  $z_2 \neq z_1$ then  $z_2(X)= \frac I X$, where $I^2=1$. Then,  $z_2$ or $z_1\, z_2$ is the
transformation $X \to \frac 1 X$; say $z_2(X)=\frac 1 X$.

Thus,  we have $\{ \pm \a_1, \pm \frac 1 {\a_1}, \dots , \pm \a_n, \pm \frac 1 {\a_n} \} \subset \P$ where
$n=\left[\frac {g+1} 2\right]$. If either $z_2$ or $z_1 z_2$ fixes two Weierstrass points then $ \pm 1$ or $
\pm I$ are also in $\P$.  Hence, the equation of $\X_g$ is given by

\begin{equation}\label{e_q}
\begin{split}
Y^2 = \left\{ \aligned
   &  \prod_{i=1}^{n} (X^4-\lambda_i X^2 +1),
 \quad   \textit{ where    } \quad    n=  \frac {g+1} 2 , \, \, g \equiv 1\mod 2    \\
   &  (X^2 \pm 1) \, \prod_{i=1}^{n} (X^4-\lambda_i X^2 +1),
 \quad   \textit{ where  } \quad    n=\frac {g} 2, \, \,g \equiv 0\mod 2    \\
   &(X^4 - 1) \, \prod_{i=1}^{n} (X^4-\lambda_i X^2 +1),
 \quad   \textit{ where} \quad    n=\frac {g-1 } 2, \, \,g \equiv 1\mod 2    \\
\endaligned \right.
\end{split}
\end{equation}
where $\l_i = \a_i^2 + \frac 1 {\a_i^2}$.  Let $s:=\l_1 + \dots + \l_n$ and  recall that $u_1:=a_1^{g+1} +
a_g^{g+1}$,   $u_g :=2 a_1 a_g$.

In the first case of the formula we have   $a_1=a_g=-s$. Then, $u_1=2s^{g+1}$ and $u_g=2 s^2$ and they
satisfy $2^{g-1}\, u_1^2=  u_g^{g+1}.$ Furthermore, they satisfy the first factor of the equation in ii). In
this case no Weierstrass points are fixed by any involutions of $V_4 \emb \bar G$, hence they lift to
involutions in $G$.

In the second case of Eq.~\eqref{e_q}, if   $X^2+1$ is a factor then $a_1=a_g=1-s$ and  $2^{g-1} u_1^2 -
u_g^{g+1}=0$.  If $X^2-1$ is a factor then
$$F(X)= X^{2g+2} - (s+1) X^{2g} + \cdots + (s+1) X^2 -1.$$
This is not in the normal form since the coefficient of $X^0$ is -1.  As in the proof of Lemma~\ref{decompo}
(cf., section 3.2) we transform the curve by
$$( X, Y )   \longrightarrow \left(  \frac 1 {(-1)^{\frac 1 {2g+2} } X},  \frac { I \cdot Y} {X^{g+1}} \right).$$
Using the formula \eqref{eq_17} (cf., section 3.2) we get
$$a_1 = \frac {s+1} { (-1)^{ \frac g {g+1} }   }, \quad a_g = - \frac {s+1} { (-1)^{\frac 1 {g+1} }    }.$$
Then,
$$u_1 = 2 (s+1)^{g+1}, \quad u_g = 2 (s+1)^2$$
and they satisfy $2^{g-1}\, u_1^2 - u_g^{g+1}$.

In the third case, one of the factors of the equation is $X^4-1$. Then, by using the same technique as above
we get
$$u_1 = -2 s^{g+1}, \quad u_g = 2 s^2$$
and the result follows. Further, they satisfy the second factor of the equation in ii). In this case each of
$z_1$ and $z_2$ fix two  Weierstrass points, hence they lift to elements of order 4 in $G$.
\end{proof}

\subsection{Computing the dihedral invariants}

The most straightforward method to decide if a hyperelliptic  curve $\X_g$ of genus $g$ defined over $k$ has
an extra involution,  and, in the affirmative case, to compute the dihedral invariants  of $\X_g$, is by
solving a polynomial system of equations.

Given  the curve $\X_g$ we want to find $\alpha = \begin{pmatrix} a & b \\ c & d \end{pmatrix} \in GL_2(k)$
such that $\X_g^\alpha$ is written in the normal form of Eq.~(2), for some $a_i \in k(a,b,c,d)$. We get a
polynomial system  by equating to zero the coefficients of $X$ to the odd power and to one the leading and
the constant  coefficients. So, we have $4$ unknowns and $g+3$ equations.  This method is simple, but
unfortunately inefficient, even for small genus  $g$.

We will present a faster method to compute the dihedral invariants if the polynomial $E(X)$ has a
decomposition. The polynomial decomposition problem can be stated as follows: given a polynomial $E \in
k[X]$, determine whether there exist polynomials $G, H $ of degree greater than one such that $E=G \circ H=G
( H(X))$, and in the affirmative case  to compute them. From the classical L\"uroth's theorem this problem is
equivalent to deciding if  there exists a proper intermediate field in the finite algebraic extension $k(E)
\subset k(X)$. From the computational point of view, there are several polynomial time algorithms for
decomposing polynomials; see \cite{Ga}. One of the main techniques  is based on the computation of the
$s$-root $H(X)$ of the polynomial $E(X)$. In that case, $\deg(E-H^s)< rs-s$, where $\deg(E) \ = \ r s$.

\begin{lemma}\label{decompo}
Let $L$ be a subfield of $k$ and  $\X_g$ a genus $g$ curve with equation $Y^2= E(X)$, where $E \in L[X]$. If the
polynomial $E(X)$ decomposes as follows:
\begin{equation}
E(X) = (G \circ H) (X), \quad  where \quad \deg (H)=2,
\end{equation}
and  $G, H\in L[X]$, then $\X_g \in \L_g$  and $\u(\X_g) \in \L_g(L)$.
\end{lemma}

\begin{proof}
If $E(X)$ has a decomposition as above, then there exist $G(X), H(X) \in L[X]$ such that $E(X)=G(H(X))$, where
$H(X)=X^2+aX$ for some $a \in L$. Let $\alpha(X)=X-a/2$ and denote by $\X_g^\alpha$ the curve after the coordinate
change $\alpha$. Then, $\X_g^\alpha$ is isomorphic to $\X_g$ and  is given by the equation
\begin{equation}Y^2=b_{g+1} X^{2g+2} + b_g X^{2g} + \dots + b_1 X^2 +b_0,
\end{equation}
 where $b_i \in L$ and $b_0\, b_{g+1} \not=0$. Without loss of generality we can assume that $b_{g+1}=1$.
 By  the following transformation $X \to b_0^{-1/(2g+2)}X^{-1} $ in $k$  the curve has the equation $Y^2=F(X)$
 where
\begin{equation}
\begin{small}
F(X)= X^{2g+2} +  {c_g} {b_0^{- \frac g {g+1}} } X^{2g} + \dots + {c_{g-i} } {  b_0^{- \frac {g-i} {g+1} } }
X^{2(g-i)} + \dots+ b_0^{- \frac 1 {g+1}} c_1 X^2 +1,
\end{small}
\end{equation}
and $c_i \in L$. The first claim follows by Lemma~\ref{lem1}. For the rest, it is straightforward to check
that the dihedral invariants of $Y^2=F(X)$ are
\begin{equation}\label{eq_17}
u_i=  \frac { c_1^{g-i+1} \, c_i} { b_0} + \frac {c_g^{g-i+1}\, c_{g-i+1}}   {b_0^{g-i+1} }
\end{equation}
for all $1 \leq i \leq g$.  Hence, $u_i\in L$.
\end{proof}

If  $E(X)$ is a tame polynomial (i.e.,  $2g+2$ is prime to the characteristic of the field $k$) then the
computation of $G(X)$ and $H(X)$ only requires $O(g^2)$ arithmetic operations in the ground field $k$; see for
instance \cite{Gu}. So, the above lemma provides an algorithm that only requires $O(g^3)$ field arithmetic
operations. If $k$ is a  zero characteristic field then a polynomial $E(X)\in F[X] $ is indecomposable over the
subfield  $F \subset k$ if and only if  $E(X)$ is indecomposable over $k$.  In particular, if the curve is defined
over the rational number field $\Q$ having  an extra involution, then the dihedral invariants are also in $\Q$.

\section{Field of moduli of curves}

In this section all curves are defined over $\C$.  For each $g$, the moduli space   $\M_g$  (resp., $\H_g$)
is the set of isomorphism classes of genus $g$ irreducible, smooth, algebraic   (resp.,  hyperelliptic)
curves $\X_g$ defined over $\C$. It is well known that $\M_g$ (resp.,  $\H_g$) is a $3g-3 $  (resp., $2g-1$)
dimensional variety. Let $L$ be a subfield of $\C$. If $\X_g$ is  a genus $g$ curve defined  over $L$, then
clearly $[\X_g]\in \M_g(L)$. Generally, the converse does not hold. In  other words,  the moduli  spaces
$\M_g$ and $\H_g$ are coarse moduli spaces.

Let $\X$ be a curve defined over $\C$.  A field $F \subset \C$  is called a   \textbf{field of definition} of
$\X$ if there exists $\X'$ defined over $F$ such that $\X'$ is isomorphic to $\X$ over $\C$.

\begin{definition}
The   \textbf{field of moduli}  of $\X$ is a subfield $F \subset \C$ such that for every automorphism $\sigma
\in \Aut (\C)$ the following holds:  $\X$ is isomorphic to $\X^\sigma$  if and only if $\, \,  \sigma_F =
id$.
\end{definition}

We will use $\p=[\X]\in   \M_g$ to denote  the corresponding   \textbf{moduli point}  and $\M_g (\p)$ the
residue field of $\p$  in $\M_g$. The field of moduli of $\X$ coincides with the residue field $\M_g (\p)$ of
the point $\p$ in $\M_g$; see Baily \cite{Ba2}.
%
The notation $\M_g (\p)$ (resp., $M(\X)$ ) will be used to denote the field of moduli of $\p \in
\M_g$ (resp., $\X$). If there is a curve $\X^\prime$ isomorphic to $\X$ and defined over $M(\X)$,
we say that $\X$ has a   \textbf{rational model over its field of  moduli}. As mentioned above, the
field of moduli of curves is not necessarily a field of definition;, see    \cite{Shi1} for
examples of such families of curves.


\subsection{Conditions for the field of moduli to be a field of definition.}
What are necessary conditions for a curve to have a rational model over its field of moduli? We consider only
curves of genus $g> 1$; curves of genus 0 and 1 are known to have a rational model over its field of moduli.
In (1954) Weil  showed that;

\medskip

{\bf i)}    \textit{For every  curve $\X$ with trivial automorphism group  the field of moduli is a field of
definition.}

\medskip

Later work of  Baily, Shimura, Coombes-Harbater, D\'ebes, Douai, Wolfart et al. has added other conditions which
briefly are summarized below.

\medskip

\noindent \textit{The field of moduli of a curve $\X$ is a field of definition if:}

\medskip

{\bf ii)}    \textit{ $Aut (\X)$ has no center and has a complement in the automorphism group of $Aut (\X)$}

\medskip

{\bf iii)}    \textit{The field of moduli $M (\X)$ is of cohomological dimension $\leq 1$ }

\medskip

{\bf iv)}   \textit{The canonical $M ( \X)$-model of $\X / Aut(\X)$ has $M (\X)$-rational points.}

\medskip

\noindent The proofs  can be found in \cite{We}.

\subsection{Field of moduli of hyperelliptic curves}

In his 1972 paper  \cite{Shi1} Shimura proved that:

\begin{theorem}[Shimura]  No generic hyperelliptic curve of even genus has a model rational over its field of
moduli.
\end{theorem}

\noindent A generic hyperelliptic curve has automorphism group of order 2. Shimura's family and Earle's
family of curves (i.e., with non-trivial obstruction) are both families of hyperelliptic curves with
automorphism group of order 2.  Consider the following problem:

\medskip

\noindent   \textbf{Problem 1:}     \textit{Let the  moduli point $\p \in \H_g$ be given. Find necessary and
sufficient conditions  such that the field of moduli $M (\p)$ is a field of  definition.  If $\p$ has a
rational model $\X_g$ over its field of moduli, then determine explicitly the equation of $\X_g$.}

\medskip

Mestre (1993) solved the above problem for genus two curves with automorphism group $\Z_2$; see \cite{Me} for
details.  Mestre's approach is followed by Cardona/Quer (2003) to prove that for points $\p \in \M_2$ such
that $| Aut (\p)| > 2$, the field of moduli is a field of definition. Algorithms have been implemented which
combine these results and give a rational model of the curve (when such a model exist) over its field of
moduli.  However, the problem is quite open for $g \geq 3$. Especially, there are no such explicit results as
in the case $g=2$. We conjecture the following:

\medskip

\noindent   \textbf{Conjecture:}  \emph{Let $\p  \in \H_g $ such  that $|Aut(\p)|  > 2$.  Then the field of
moduli of $\p$ is a field of definition.}

\medskip

\begin{remark} The above was first conjectured during a talk of the second author in ANTS V (Sydney, 2002).
For the first time in print it has appeared in \cite{Sh3}.
\end{remark}

In studying the above conjecture it becomes important to first determine a list of groups that
occur as automorphism groups of genus $g$ curves for a given $g$. The most up to date work on this
is \cite{MS} where explicit lists are provided for small $g$. The automorphism groups of
hyperelliptic curves have been studied in \cite{BS, Bu}. For a complete list of such groups and
algorithms computing the automorphism group of a given curve see \cite{Sh4}. Next we prove the
conjecture for all moduli points $\p \in \L_g$ such that $V_4 \emb \bAut (\p)$.

\begin{theorem}
If $\p= (u_1, \dots , u_g) \in \L_g$  such that $ 2^{g-1}\, u_1^2 - u_g^{g+1} =0 $
then the field of moduli is a field of definition. Moreover, the rational model over the field of moduli is
given by
\begin{equation}\label{rat_mod}
\X_g: \quad Y^2=u_1 X^{2g+2} + u_1 X^{2g} + u_2 X^{2g-2} + \dots +  u_g X^2 +2.
\end{equation}
\end{theorem}

\begin{proof}
Let $\p=(u_1, \dots , u_g)\in \L_g$ such that $V_4 \emb \bAut(\p)$. Hence, $2^{g-1} u_1^2=u_g^{g+1}$.  All we
need to show is that the curve $\X_g$ given in Eq.~\eqref{rat_mod} corresponds to the moduli point $\p$. By
an appropriate transformation $\X_g$ can be written as
\begin{equation}Y^2= X^{2g+2} + ( \frac {u_1} 2 )^{\frac 1 {g+1}} \cdot X^{2g} +
\sum_{i=1}^{g-1} \frac {u_{g+1-i}} {u_1} \cdot (\frac {u_1} 2 )^{\frac {g+1-i} {g+1}}  \cdot X^{2i}+1.
\end{equation}
Then, its dihedral invariants are
\begin{equation}
\begin{split}
u_1(\X_g) &= \frac {u_1} 2 + (\frac {u_g} {u_1} )^{g+1} \cdot (\frac {u_1}2)^g= \frac { 2^{g-1} u_1^2+u_g^{g+1}
} {2^g u_1}, \\
u_j (\X_g) & = u_j, \quad \textit { for  } \quad j=2, \dots ,\,  g.\\
\end{split}
\end{equation}
Substituting $u_g^{g+1}=2^{g-1} u_1^2$ we get  $u_1(\X_g)=u_1$. Thus, $\X_g$ is in the  isomorphism class
determined by $\p$.   Because   coefficients of $\X_g$ are given as rational functions of $u_1, \dots , u_g$
the curve is defined over its field of moduli. This completes the proof.
\end{proof}

\begin{cor}
Let $\p \in \H_g$ such that $V_4 \emb \Aut(\p)$. Then the field of moduli of $\p$ is a field of definition
with rational model as in Eq.~\eqref{rat_mod}.
\end{cor}

We illustrate next with cases $g=2, 3$.  The case $g=2$ is the only case which is fully understood.
\begin{lemma}
Let $ \u  \in \M_2  $ such that $|\Aut (\u) | > 2$.  Then, the field of moduli of $\u$ is a field of definition.
Moreover, a rational model over the field of moduli is given by:

\medskip

\subitem i) If $Aut(\u)\iso D_8$ then
$$Y^2= u_1 X^6 + u_1 X^4 +  u_2 X^2 +2.$$

\subitem ii) If $Aut(\u)\iso D_{12}$ then
$$Y^2=    4(u_2-450)\, X^6 + 4 (u-2-450) \, X^3 + u_2-18$$

\subitem iii) $Aut(\u)\iso V_4$

\subsubitem a)   If  $u_2\neq 0$ then
\begin{small}
\begin{equation*}
\begin{split}
Y^2  = & \, \, \frac 8 {d_6^3}  (u_2^3+u_2^2u_1 + 2d_6) X^6  \,
+ \,  \frac 8 {d_6^2} (u_2^2+12 u_1) X^5 \,         \\
 & + \, \frac 4 {d_6^2} (15 u_2^3 - u_2^2 u_1 + 30 d_6) x^4 \,
- \, \frac 8 {d_6} (u_2^2  - 20 u_1) X^3 \,         \\
 & + \, \frac 2 {d_6^2} (15 u_2^3 - u_2^2 u_1 + 30 d_6) X^2 \,
+ \, 2 (u_2^2 + 12) X \, + \, (u_2^3+u_2^2u_1 + 2d_6)
\end{split}
\end{equation*}
\end{small}
where $d_6=2u_1^2-u_2^3$.

\subsubitem b)   If  $u_2 = 0$ then
\begin{small}
\begin{equation*}
\begin{split}
Y^2  = & \, \, (2u_1+1)X^6 - 2(4u_1-3)X^5+(14u_1+15)X^4 - 4(4u_1 - 5)X^3\\
 &  +(14u_1+15)X^2 - 2(4u_1-3)X + 2u_1+1.
\end{split}
\end{equation*}
\end{small}
\end{lemma}

\proof  For parts i) and ii) see \cite{Sh2}. For iii) compute the absolute invariants $i_1, i_2, i_3$ and
check that they are the same as in expressions in equation (19) in \cite{Sh-V}. Hence, the dihedral
invariants are $u_1, u_2$ since they provide a birational parametrization of the space $\L_2$.

\qed

Part i) and ii) of the Lemma were proved in \cite{Sh2}. Part iii) was the main focus of \cite{CQ}. The
approach there, however, uses absolute invariants and  the equation of the curve is more complicated.
The reader should compare the equations of the above lemma with those provided in \cite{CQ} in order to be
convinced of the advantages of using the dihedral invariants. For  $g=3$,   we have the following.

\begin{lemma} Let $ \u  \in \L_3 (k) $ such that $|Aut( \u )| > 4$.
Then, there exists a genus 3 hyperelliptic curve $\X_3$ defined over $k$ such that $\u(\X_3)=\u$. Moreover,
the equation of $\X_3$ over its field of moduli is given by:

\medskip

\subitem i) If $|Aut(\X_3)|=16$ then $$Y^2= w X^8 +w X^4 +1.$$

\smallskip

\subitem ii) If $Aut(\X_3)\iso  D_{12} $ then
\begin{small}
\begin{equation*}
\begin{split}
Y^2 = & \, \, (u_3-260) X^8 - 7 (u_3-98)X^6 + 15 (u_3-134)X^4 \\
& - 9(u_3-162)X^2 + 126
\end{split}
\end{equation*}
\end{small}
where $u_1, u_2, u_3$ satisfy equations (14).

\smallskip

\subitem iii) If $Aut\iso \Z_2 \times \Z_4$ then
$$Y^2=u_3^4 X^8+u_3^4X^6+8 u_3X^2-16.$$

\smallskip

\subitem iv) If $Aut(\u)\iso \Z_2^3$ then
$$ Y^2= u_1 X^8 +u_1 X^6 + u_2 X^4 +  u_3 X^2 +2. $$
\end{lemma}

\proof The proof in all cases consists of simply computing the dihedral invariants. It is easy to check that
these dihedral invariants satisfy the corresponding relations for $Aut(\X_3)$ given in \cite{Sh9}.

\qed

\begin{cor}
Let $\p \in \H_3$ such that $| \Aut (\p) | > 4$. Then the field of moduli of $\p$ is a field of definition.
\end{cor}

\proof There is only one hyperelliptic curve of genus 3 which has no extra involutions and order of the
automorphism group $> 4$; see \cite{MS} or \cite{Sh9}.  This curve is $Y^2= X^7-1$ and its field of moduli is
$\Q$. The result follows from the above Lemma.

\qed

\section{Concluding Remarks}

The main goal of this paper was to introduce  dihedral invariants and  show how they can be used to answer
some classical problems. In \cite{Sh4} we use such invariants to design an algorithm which determines the
automorphism group of hyperelliptic curves. In section four we give another example of such applications.

The field of moduli problem discussed in section four is a classical problem of algebraic geometry. There are
many works in the literature which extend the problem to other categories other than curves (i.e., covers,
polarized abelian varieties, etc.). However, none of these papers gives an explicit way of determining the
field of moduli or providing a rational model of the curve over the field of moduli when such a model exist.
Dihedral invariants are useful in this direction when dealing with hyperelliptic curves with extra
involutions.

For $g > 3$ providing rational models over the field of moduli is a difficult task. In \cite{Sh6} such models
are provided for all hyperelliptic curves $\X_g$ of genus $g \leq 12$ and $\bAut (\X_g) \iso A_4$.



\begin{thebibliography}{00}

\bibitem {Ba}  W. Baily, On the automorphism group of a generic curve of genus $>2$,
{\it J.  Math. Kyoto Univ.} {\bf 1} (1961/1962), 101--108; correction, 325.


\bibitem  {Ba2}    W. L. Baily, On the theory of theta functions, the moduli of abelian varieties and
the moduli of curves, Ann. Math. 75 (1967), 342--381.


\bibitem{BS}   R. Brandt and H. Stichtenoth, Die Automorphismengruppen hyperelliptischer Kurven,
{\it Man. Math} {\bf 55} (1986), 83--92.


\bibitem{Bu} E. Bujulance, J. Gamboa, and G. Gromadzki,
The full automorphism groups of hyperelliptic Riemann surfaces, {\em Manuscripta Math.} {\bf 79} (1993), no.
3-4, 267--282.

\bibitem {CQ}   G. Cardona and J. Quer, Field of moduli and field of definition for curves of genus 2,
Article math.NT/0207015








\bibitem  {DM}   P. Deligne P, D. Mumford D,
The irreducibility of the space of curves of given  genus, Publ. Math.Hautes \'Etudes Sci. {\bf 36}, 75-109, 1969.



\bibitem {Ga}  J. von zur Gathen, Functional decomposition of polynomials: the tame case, J. of Symbolic
Computation, 9 (1990), 281- 299.

\bibitem {Gu}   J. Gutierrez, A polynomial decomposition algorithm over factorial domains,
Comptes Rendues Mathematiques, de Ac. de Sciences, 13 (1991), 81-86.

\bibitem {vishi} V. Krishnamorthy, T. Shaska, and H. V\"olklein;
Invariants of binary forms, \emph{Developments in Mathematics}, Kluwer Academic Publishers, Boston, MA, 2004,
101 -122.


\bibitem {MS}   K. Magaard, T. Shaska, S. Shpectorov, and H. V\"olklein,  The locus of curves with prescribed
automorphism group. Communications in arithmetic fundamental groups (Kyoto, 1999/2001).
S\=urikaisekikenky\=usho K\=oky\=uroku No. 1267 (2002), 112--141.


\bibitem {Ig} {\sc J. Igusa}, Arithmetic Variety of Moduli for genus 2. {\it Ann. of Math}. (2), 72, 612-649,
(1960).


\bibitem {Me}   P. Mestre,   Construction de courbes de genre 2 \'a partir de leurs modules. In T.   Mora and
C. Traverso, editors,  {\it Effective methods in algebraic geometry}, volume 94. {\it Prog.  Math.}, 313-334.
Birkh\"auser, 1991. Proc. Congress in Livorno, Italy, April 17-21, (1990).









\bibitem {Sh2} T. Shaska,
Genus 2 curves with (3,3)-split Jacobian and large automorphism group, Algorithmic Number Theory (Sydney,
2002), {\bf 6}, 205-218, Lect. Not. in Comp. Sci., 2369, Springer, Berlin, 2002.


\bibitem {Sh3} T. Shaska,
Computational aspects of hyperelliptic curves, Computer mathematics. Proceedings of the sixth Asian symposium
(ASCM 2003), Beijing, China, April 17-19, 2003. River Edge, NJ: World Scientific. Lect. Notes Ser. Comput.
10, 248-257 (2003).


\bibitem {Sh4} T. Shaska,
Determining the automorphism group of hyperelliptic curves, {\it Proceedings of the 2003
International Symposium on Symbolic and Algebraic Computation}, ACM Press, pg. 248 - 254, 2003.




\bibitem {Sh6} T. Shaska,
Some special families of hyperelliptic curves, {\it J. Algebra Appl.}, vol \textbf{3}, No. 1 (2004), 75-89.



\bibitem {Sh-V} T. Shaska, H. V\"olklein;
Elliptic subfields and automorphisms of genus two fields, {\it Algebra, Arithmetic and Geometry with
Applications. Papers from Shreeram S. Abhyankar's 70th Birthday Conference}, pg. 687 - 707, Springer (2004).



\bibitem {Sh10} T. Shaska and J. L. Thompson,
On the generic curves of genus 3, {\it Affine Algebraic Geometry}, Cont. Math., American
Mathematical Society, 2004.

\bibitem {Sh9} T. Shaska,   Hyperelliptic curves of genus 3 with split Jacobians, (work in progress).

\bibitem {Shi1}   G. Shimura,
On the field of rationality for an abelian variety. Nagoya Math. J. 45 (1972), 167--178.



\bibitem  {V} {\sc H. V\"olklein}, \emph{Groups as Galois groups. An introduction.}  Cambridge Studies in Advanced
Mathematics, 53. Cambridge University Press, Cambridge, 1996.



\bibitem  {We}  A. Weil, The field of definition of a variety, Amer. J. Math. 78 (1956), 509--524.


\end{thebibliography}
\end{document}